\documentclass[11pt]{amsart}

\usepackage{amsmath}
\usepackage{amssymb}
\usepackage{graphicx}
\usepackage{amscd}
\usepackage{stmaryrd}
\usepackage[mathscr]{euscript}

\theoremstyle{definition}

\numberwithin{equation}{section}

\newcommand{\sB}{\mathscr B}
\newcommand{\sC}{\mathscr C}

\newcommand{\sH}{\mathscr H}
\newcommand{\sJ}{\mathscr J}
\newcommand{\sK}{\mathscr K}

\newcommand{\sR}{\mathscr R}

\title[Update on the Quasicentral Modulus
]{Perturbations of Operators and Non-Commutative Condensers, an Update on the Quasicentral Modulus  }

\begin{document}
\author[D.-V. Voiculescu]{Dan-Virgil Voiculescu${}$}

\address{Department of Mathematics \\ University of California at Berkeley \\ Berkeley, CA\ \ 94720-3840}
\email{{\tt dvv@math.berkeley.edu}}




\begin{abstract}
This is an update on the quasicentral modulus, an invariant for an n-tuple of Hilbert space operators and a rearrangement invariant norm, that plays a key-role in sharp multivariable generalizations of the classical Weyl - von Neumann - Kuroda and Kato - Rosenblum theorems of perturbation theory. There are also connections with self-similar measures on certain fractals and to the Kolmogorov-Sinai dynamical entropy. Some open problems are also pointed out. Recently a non-commutative analogy with condenser capacity in nonlinear potential theory is emerging,that  provides a new perspective on the subject.
\end{abstract}

\maketitle


\section{Introduction}
\label{sec1}
A perturbation of an operator, or more generally of an n-tuple of operators, is usually assumed to be small in some sense. There is a big difference about what "small" can mean, between the finite dimensional case and the infinite dimensional case. In finite dimension there is only one way for a perturbation to be small, while in infinite dimension the existence of non-trivial ideals of the ring of bounded operators on a Hilbert space provides many different possibilities. In finite dimension a small perturbation of a selfadjoint operator produces small changes in the eigenvalues, while in infinite dimension there are new phenomena. It is a consequence of the Kato-Rosenblum theorem that a trace-class perturbation preserves up to unitary equivalence the Lebesgue absolutely continuous part of the selfadjoint operator and the Weyl-von Neumann -Kuroda theorem implies that the result is sharp, in that for larger ideals than the trace-class there are perturbations that destroy the Lebesgue absolutely continuous part.

\bigskip

A key ingredient in our approach to sharp multivariable generalizations of the one operator perturbation facts is the quasicentral modulus, a numerical invariant associated with an  $n$-tuple  of bounded Hilbert space operators and a normed ideal of compact operators. It has turned out the invariant has also other interesting properties, beyond perturbations, for instance we found that for an appropriate ideal there are connections to the Kolmogorov-Sinai dynamical entropy.

\bigskip

Wondering about the nice properties of the quasicentral modulus, we recently found that it has an extension that can be viewed as a non-commutative nonlinear capacity, that is that some kind of noncommutative nonlinear potential theory may be at work, underlying the perturbation theory results.

\bigskip

This survey proceeds as follows. We begin with some of our early sharp multivariable perturbation results and point out some important open problems for perturbations in dimension two, that is for pairs of commuting selfadjoint operators. Next we go over to certain results about the quasicentral modulus. In the case of the so-called Macaev ideal, which is related by duality to the ideal where the Dixmier trace lives, we have the connection with the Kolmogorov-Sinai entropy.  Then in the case of $n$-tuples of commuting selfadjoint operators there is an exact  formula for the quasicentral modulus with respect to the corresponding threshold ideal for perturbations . This formula has had applications in noncommutative geometry. We then pass to recent results on exact formulae: the exact formula for the hybrid quasicentral modulus , where the normed ideal is replaced by an $n$-tuple of normed ideals, that is the invariant for perturbations of $n$-tuples of operators, where the perturbation of each component is from  a different ideal. Another situation considered recently, is the exact formula for commuting $n$-tuples with joint spectrum in a certain fractal. 

\bigskip

Then we turn to the nonlinear potential theory noncommutative analogy. We first recall the definition of capacity for condensers in the nonlinear Sobolev setting which we will use. Then we provide a dictionary of the noncommutative analogy and  we point out that in order to improve the analogy it is necessary to pass to the nonlinear potential theory on graphs. In the case of Cayley graphs of finitely generated groups we have actually that capacity quantities for certain condensers actually coincide in the analogy. A natural next step to explore the analogy is to take a look at noncommutative variational problems associated with noncommutative condensers. What we found in the case of the Schatten-von Neumann  $p$-class, points at a noncommutative analogue of the $p$-Laplace equation. 

\bigskip

Another natural question is how far the perturbation theory can be extended beyond the type $I_{\infty}$ case of $\sB(\sH)$ to the semifinite setting which includes type $II_{\infty}$ factors. At this time, it is not clear whether the complexity of the type $II$ von Neumann algebra world will play a role in this, so we only mention a few papers where the partial results available at present can be found.

\bigskip

An aspect  of the quasicentral modulus, not covered in this survey, is the relation to commutants mod normed ideals which we surveyed in \cite{25}. Note that commutants mod normed ideals are a certain class of noncommutative smoothness spaces and that the appearance of smoothness spaces along with some potential theory, noncommutative in this case, should not be surprising.

\bigskip
\noindent

\section{Preliminaries on normed ideals }
\label{sec2}

Let  $\sH$ be a separable infinite dimensional complex Hilbert space and let $\sB(\sH) $ be the von Neumann algebra of bounded operators on $\sH$. By $\sR(\sH)$, $\sK(\sH)$ or simply $\sR$ and $\sK$ we will denote the ideals of finite rank and respectively of compact operators on $\sH$. Also $\sR^+_1 (\sH)$ or $\sR^+_1$ will stand for $ \{ A \in \sR \vert 0 \leq A \leq I \} $ the finite rank positive contractions. By $\sB(\sH)_{sa}$, $\sK_{sa}$, $\sR_{sa}$ etc. we will denote the corresponding subsets of selfadjoint elements.

\bigskip

Roughly a normed ideal $(\sJ, \vert \cdot \vert_{\sJ})$ can be described as an ideal of $ \sB (\sH)$, $\sJ \subset \sK$, which is a Banach space with respect to a norm $\vert \cdot \vert _{\sJ}$, which is of the form 

\[
{\vert T \vert_{\sJ} = \vert (s_j)_ { j \in \mathbb {N}} \vert _{\Phi }}
\]

where 

\[
{s_1 \geq s_2 \geq ... }
\]
\noindent
are the eigenvalues of $(T^*T)^{1/2}$ and $\vert \cdot \vert_{\Phi} $ is a rearrangement invariant norm on sequences. (For more precise details see \cite {10} or \cite {16} or \cite {15} which covers the more general semifinite generalizations of these concepts.)

\bigskip
\noindent

Note that in particular if $A, B \in \sB(\sH) $ and $X \in \sJ$ we have

\[
{\vert AXB \vert _{\sJ} \leq ||A|| \vert X \vert _{\sJ } ||B|| }.
\]

\bigskip
\noindent
If $1 \leq p < \infty$ and 

\[
{ \vert (s_j)_{j\in \mathbb {N}}\vert_{\Phi}  = \vert (s_j)_{j\in \mathbb {N}}\vert_p = \left( \sum_{j}  s^p_j \right) ^{1/p}}
\]

we have the Schatten-von Neumann $p$-class $(\sC_{p}, \vert \cdot \vert _p) $ . For us here the Lorentz $(p, 1)$ ideals given by 

\[
{\vert (s_j)_{j\in \mathbb {N} }\vert_{\Phi}  =  \sum_{j} j^{-1 + 1/p} s_j } 
\]

will be particularly important. In case $p = \infty$ , this gives

\[
{\vert (s_j)_{j\in \mathbb {N} } \vert_{\Phi} = \sum_{j} j^{-1} s_j} .
\]

We will denote these normed ideals by $ (\sC^-_p, \vert \cdot \vert^-_p)$ and $(\sC^-_{\infty} , \vert \cdot \vert ^-_{\infty} )$ will be called the Macaev ideal following \cite{10}.
These ideals are also refered to as the Lorentz $(p,1)$ ideals and we have $\sC ^-_p \subset \sC_p $ and $\sC^-_1 = \sC_1$ .

\section {The quasicentral modulus $k_{\sJ}(\tau)$ }
\label{sec3}  
\bigskip
\noindent

If $(\sJ, \vert \cdot \vert_{\sJ}) $ is a normed ideal and $ \tau = (T_j ) _{1 \leq j \leq n} \in (\sB (\sH))^{n} $ is an $n$-tuple of bounded operators, the quasicentral modulus is defined as follows:\\
\\
\[
{k_{\sJ}(\tau) = \text{ least } C \in [ 0, \infty] \  \text{so that}  \exists A_m \in \sR^{+}_{1}, m\in \mathbb {N}, A_m  \uparrow I} 
\]   

and
\[
{\vert [A_m , \tau ] \vert_{\sJ} \rightarrow C \text{ as } m \to  \infty }
\]

where
\[
{\vert [ A_m , \tau ] \vert _{\sJ} = \max_{ 1 \leq j \leq n}  \vert [A_m , T_j] \vert _{\sJ} }.
\]

\bigskip

A basic property of $k_{\sJ} (\tau) $ is that if $\sR$ is dense in $\sJ$ in $\sJ$ - norm , then 
$k_{\sJ} (\tau) = k_{\sJ} (\tau ' ) $ if $ \tau \equiv \tau ' $ mod $\sJ$ .

\section{The sharp analogue for $n$-tuples $n \geq 3 $ of the Kato-Rosenblum theorem }
\label{sec4}

\bigskip
\noindent
{\bf Theorem}
\cite {19} Let $\tau , \tau'$ be two $n$-tuples of commuting bounded selfadjoint operators, $n \geq 3$, so that $\tau \equiv \tau'$ mod $\sC^-_{n}$. Let $\sH_{ac} (\tau), \sH_{ac}(\tau')$ be the Lebesgue absolutely continuous subspaces of $\tau$ and $\tau'$ and let $\tau_{ac} = \tau \vert \sH_{ac} (\tau), \tau' \vert \sH_{ac}(\tau') $. If $u_m \in \sC^{\infty} (\mathbb {R} ^n), m \in \mathbb {N} $ are so that $|u_m| = 1$ in a compact neighborhood $K$ of $\sigma (\tau) \cup \sigma (\tau')$ and $w-\lim_{m\rightarrow \infty} u_m = 0$ in $ L^2 (K, d\lambda) $, $\lambda$ denoting Lebesgue measure, then the strong limit 

\[
{W = s-lim_{m\rightarrow \infty} u_m (\tau')^{\ast} u_m (\tau) \vert \sH_{ac} (\tau) }
\]

exists and is a unitary intertwiner of $\tau_{ac} $ and $\tau'_{ac} $.

\bigskip

{\bf Remarks} : a) $W$ does not depend on the choice of the sequence $(u_m)_{m \in \mathbb {N} }$. Indeed two different sequences can be mixed and form a sequence to which the theorem applies, hence the $W$  for the two sequences are equal.

 \bigskip
 
 b) If $\tau$ and $\tau'$ are $n$-tuples of commuting bounded selfadjoint operators so that $\tau \equiv \tau'$ mod $ \sC^{-}_n$, $n \geq 3$ then $\tau_{ac} $ and $\tau'_{ac} $ are unitarily equivalent.
 
 \bigskip
  
 c) In \cite{19} there is also a weaker version of the the theorem for $n \geq 2$ which still implies that $\tau_{ac} $ and $\tau'_{ac} $ are unitarily equivalent, however the limit defining $W$ is not a strong limit. We discuss in section 5 the open problems which arise. 

\bigskip

d) The reason the Theorem and Remark b) are sharp statements is due to another result obtained using the quasicentral modulus \cite {3}, a sharp multivariable analogue of the Weyl-von Neumann -Kuroda theorem : if $\sJ $ is a normed ideal $\sJ \not \subset \sC^{-}_n$, then given an $n$-tuple of commuting selfadjoint operators $\tau$ , there is another such $n$-tuple $\tau'$ so that $\tau \equiv \tau' $ mod $\sJ$ and $\tau'$ has an orthonormal basis of eigenvectors.

\bigskip                                                                                               

e) Preceding our work based on the quasicentral modulus a weaker analogue of the Kato-Rosenblum theorem for $n$-tuples, $n \geq 3$ and perturbations in some $\sC_p$ , $p < n $ had been obtained in \cite{29}. This result is not sharp and note that also here the case $n = 2$ is not settled.

\bigskip 

 f) Typical functions $u_m$  to which the Theorem applies are $u_m (x) = \exp (i \langle x , \xi^{(m)} \rangle) $ where $\xi^{(m)} \in \mathbb {R} ^{n}$ are so that $ |\xi^{(m)} | \to \infty$ as
 $m \to \infty$.

\bigskip

g) When dealing with unbounded selfadjoint operators various ways to pass to bounded operators are known and these can also be applied to $n$-tuples of commuting selfadjoint operators. This can also be used to apply the Theorem to unbounded operators.

\section { $k^{-}_{\infty} $ and entropy (\cite{21}, \cite{22} )}
\label {sec 5}

\bigskip
\noindent

Let $\theta $ be an ergodic measure-preserving automorphism of a probability measure space $(\Omega , \Sigma , \mu), \mu (\Omega) = 1 $. Let $h(\theta)$ be the Kolmogorov-Sinai entropy of $\theta$. We will denote by $U_{\theta}$ the unitary operator in $L^{2}(\Omega, \Sigma, \mu) $ induced by $\theta$. If $\Phi$ denotes the multiplication operators in $L^{2} (\Omega, \Sigma, \mu)$ by real-valued measurable functions with finite range, we define 

\[
{\sH_{P} (\theta) = \sup_{\phi \subset \Phi, \phi finite} ( k^{-}_{\infty} ( \phi \cup \{ U_{\theta} \} ) }
\]

the  { \bf perturbation entropy of $\theta$}.

\bigskip
\noindent
{\bf Theorem}. There are universal constants $0 < C_1 < C_2 < \infty $ and $ 0 < \gamma < \infty $ so that 

\bigskip

$(i) C_{1} H_{P} (\theta) \leq h(\theta) \leq C_{2} H_{P} (\theta) $ 

\bigskip

$(ii)$ if $\theta$ is a Bernoulli shift then  $\sH_{P} (\theta) = \gamma h(\theta) $ .

\bigskip

{\bf Remarks} . a) It is an open problem whether $h(\theta) = \gamma \sH_{P} (\theta) $ for all $\theta$ (i.e. not only for Bernoulli shifts).

\bigskip

b) The results about $\sH_{P} (\theta)$ were obtained using the full power of the mathematical machinery around the Kolmogorov-Sinai entropy. It would be interesting whether a more operator-theoretic approach to these results can be found.

\section {An exact formula for $k_{n} (\tau)$ in the case of commuting $n$-tuples of selfadjoint operators}.
\label{sec 6}

\bigskip
\noindent
{\bf Theorem}
\cite{18}  If $\tau$ is an $n$-tuple of commuting bounded selfadjoint operators, there is a universal constant $\gamma_{n} \in (0, \infty) $ so that

\begin{center}
\[
(k^{-}_{n} (\tau))^{n} = \gamma _{n} \int_{\mathbb R ^{n}} m(s) d\lambda(s) 
\]
\end{center}

where $m$ is the multiplicity function of $\tau$ and $\lambda$ is Lebesgue measure. If $n = 1$ , $\gamma_{1} = 1/ \pi $. 

\bigskip

{\bf Remark.}  In \cite{5} the formula was used in the characterization of classical compact differentiable manifolds as noncommutative manifolds in the sense of noncommutative geometry.

\section { The hybrid generalization of the quasicentral modulus and of the exact formula}
\label{sec 7}
\bigskip
\noindent

Let $(\sJ_{j},  | \cdot |_{\sJ_{j}} ) , 1 \leq j \leq n)$ be normed ideals and let $\tau = (T_{j} )_{1 \leq j \leq n}$ be an $n$-tuple of bounded operators. The {\bf hybrid quasicentral modulus } is defined as follows:

\[
k_{\sJ_{1} , ... , \sJ_{n}}(\tau) = \text{ least } C \in [ 0, \infty]  \text{ so that } \exists A_m \in \sR^{+}_{1}, m\in \mathbb {N}, A_m  \uparrow I    
\]
and
\[
\max_{1 \leq j \leq n} \vert [A_m , T_{j} ] \vert_{\sJ_{j}} \rightarrow C \text{ as }m \to  \infty .
\]

\bigskip

If $\sJ_{j} = \sC^{-}_{p_{j}} , 1 \leq j \leq n$ we will use the notation $k^{-}_{p_{1} , ... , p_{n}} (\tau)$ .

\bigskip
\noindent
{\bf Theorem}
( \cite{23} ) . If $\tau $ is an $n$-tuple of commuting bounded selfadjoint operators, $ 1<p_j , 1 \leq j \leq n , p^{-1}_{1} + ... + p^{-1}_{n} = 1 $ and $m$ is the multiplicity function of $\tau$ we have

\[
(k^{-}_{p_{1} , .... , p_{n} } (\tau) )^{n} = \gamma_{p_{1} , ... , p_{n} } \int _{\mathbb R^{n}} m(s) d\lambda (s) 
\]

where $ 0 < \gamma_{p_{1} , ... , p_{n} } <  \infty $ is a universal constant.

\bigskip

There are also other kind of results in the hybrid setting which can be found in \cite{24}.

\section{The exact formula in the case when the spectrum is in  certain self-similar fractals}
\label{sec 8} 
\bigskip
\noindent

The exact formula for $k^{-}_{n} (\tau) $ where $\tau$ is an $n$-tuple of commuting selfadjoint operators arises from viewing the spectrum $\sigma (\tau ) $ as a subset of $\mathbb {R}^{n} $. In case $\sigma (\tau) $ is contained in certain self-similar fractals $X$ there is an analogous exact formula:

\begin{center}
\[
(k^{-}_{p} (\tau))^{p} = \gamma_{X} \int_{X} m(s) d\sH_{P} (s) 
\]
\end{center}

where the reference measure on $X$ is the Hutchinson ( which is Hausdorff $p$-measure with $p$ the Hausdorff dimension in the case of the fractals for which such an exact formula has been proved). We refer to the books by Falconer for background on fractals.

\bigskip

In \cite{26} we proved this kind of exact formula for certain Cantor-like self-similar fractals $X$. Then Ikeda and Izumi in \cite{11} went further and proved that the exact formula also holds for a more general class of self-similar fractals which includes the Sierpinski gasket and the Sierpinski carpet. In \cite{11} there are also results for other fractals with the exact formula replaced by inequalities. Very recently Glickfield in \cite{9} showed the exact formula also holds for other fractals assuming the spectral measure of the $n$-tuple satisfies an absolute continuity condition.

\bigskip

{\bf Remark.}  We did not discuss analysis aspects underlying the quasicentral modulus results for commuting $n$-tuples of commuting selfadjoint operators. The reader may get some idea about these from the papers \cite{8} , \cite{20}.

\section{Some open problems when $n = 2$}
\label{sec9} 
 \bigskip
 \noindent
 
 The general operator algebra theorems in (\cite{18} , \cite{19}) on the existence of wave operators underlie the Theorem and Remarks in section 4 and have also many other corollaries when combined with other results on the quasicentral modulus. There are however also questions which possibly will require work from scratch, which may be the case in particular
with the question whether the results for $n \geq 3$ also hold for $n = 2$. We present here two problems in this vein. For the sake of simplicity the wave operators $W$ will be 
strong limits of exponentials. 

\bigskip

{\bf Problem 1.} Let $\alpha = (A_1 , A_2), \beta = (B_1 , B_2)$  be two pairs of commuting bounded selfadjoint operators so that $A_j - B_j \in \sC^{-}_2 , j = 1, 2. $ Assume that the spectral measure of $\alpha$ is Lebesgue absolutely continuous. Does it follow that the strong limit 

\[
{ s-\lim_{t \rightarrow + \infty} \exp(-itB_1) \exp(itA_1) }
\]

exists ?

\bigskip

{\bf Problem 2.} Assume that $\alpha$ and $\beta$ satisfy the assumptions in Problem 1. and that for some $0 \leq \theta_1 \neq \theta_2 < 2\pi $ the strong limits

\[
{ W_k = s-\lim_{t \rightarrow + \infty} \exp(-it((\cos \theta_k) B_1 + (\sin \theta_k) B_2)) \exp(it((\cos \theta_k) A_1 + (\sin \theta_k) A_2))}
\]

exist for k = 1, 2. Does it follow that $W_1 = W_2$ ?.

\bigskip 

 {\bf Remark } In the case of negative answers explicit examples will be of interest.
 
 \section{Monsieur Jourdain realization}
 \label{sec 10 } 

Why does the quasicentral modulus have nice properties, does it resemble some other interesting quantities ? Recently, perhaps like Moliere's Monsieur Jourdain, after many years I realized that the quasicentral modulus may be a noncommutative relative of the condenser capacity in nonlinear potential theory. With the quasicentral modulus, I  may have had for a long time the benefits from a bit of noncommutative nonlinear potential theory , without being aware of it.

\section{Condenser capacity in nonlinear potential theory}
\label{sec 11}

We recall here definitions about condenser capacity in nonlinear potential theory in the form convenient for the noncommutative analogy.

\bigskip

Let $\Omega \subset \mathbb{R}^{n} $ be a non-empty open set, a condenser will be a pair of disjoint compact subsets $K , L  \Subset \Omega , K \cap L = \emptyset $. The $p$-capacity $1 \leq p < \infty$ is defined by
\begin{center}
\[
 cap_p (K, L ; \Omega) = \inf \{|| \nabla u ||^{p}_{p} \vert u \in {\sC}^{\infty}_{0} (\Omega), 0 \leq u \leq 1, u \vert _{K} \equiv 1 , u \vert _{L} \equiv 0 \}  
\]
\end{center}
where $\nabla$ is the gradient, $|| \cdot ||_p$ the $p$ - norm and ${\sC}^{\infty}_{0} ( \cdot) $ denotes the infinitely differentiable functions with compact support.This also involves choosing a norm on $\mathbb{R}^{n}$ since the gradient takes values in $\mathbb{R}^{n}$ and different norms give different but comparable capacities.

\bigskip

If $L = \emptyset$ we get the $p$-capacity of $K \Subset \Omega$

\[
{cap_p (K ; \Omega) = cap_p (K, \emptyset ; \Omega)}.
\]

The $p$-capacity of the set $\Omega$ is then 

\[
{cap_p (\Omega) = \sup_{K \Subset \Omega} cap_p (K; \Omega)}.
\]

\bigskip

If $p = 2$ and the norm on $\mathbb{R}^{n}$ is the Euclidean norm, then $|| \nabla u ||^{2}_{2}$ is the Dirichlet integral and we are in the realm of linear potential theory , the starting point for the later developments. Note that linear potential theory has a farreaching noncommutative extension (see \cite{4}).

\bigskip

Replacing the $p$-norm $|| \cdot || $ by the $(p,q)$-norm $|| \cdot ||$ ( see the papers by Costea and Maz'ya \cite{6}, \cite{7} and references therein)  also suggests using general rearrangement invariant norms. Proceeding in this direction note that a power-scaling changes the defined condenser capacity, but without a loss of information since performing a reverse power-scaling will ondo the first power-scaling. The power-scaled quantity

\[
{\inf \{|| \nabla u||_p \vert  u \in {\sC}^{\infty}_{0} (\Omega), 0 \leq u \leq 1, u \vert _{K} \equiv 1, u \vert_{L} \equiv 0 \} }
\]

immediately leads to consider

\[
{\inf \{ || \nabla u ||_{\Phi} \vert u \in {\sC}^{\infty}_{0} (\Omega), 0 \leq u \leq 1, u \vert_{K} \equiv 1 , u \vert_{L} \equiv 0 \}}
\]
where $ || \cdot ||_{\Phi} $ is any rearrangement invariant norm.

\section{The dictionary for the noncommutative analogy (\cite{27}) }
\label{sec 12}

1. If $K , L  \Subset \Omega , K \cap L = \emptyset $ is a condenser , a noncommutative analogue will be a pair of orthogonal finite rank projections $P, Q , PQ = 0$.

\bigskip

2. The noncommutative analogue of the gradient $\nabla$ will be the "gradient of inner derivations " $[ \cdot , \tau] $ where $\tau = (T_j)_{1 \leq j \leq n} $ is an $n$-tuple of bounded operators. 

\bigskip

3. The correspondent of ${\sC}^{\infty}_{0} (\Omega)$ in the noncommutative analogy will be $\sR$, the finite rank operators.

\bigskip

4. The rearrangement invariant norm $|| \cdot ||_{\Phi} $ will have as noncommutative analogue, not so surprising, the norm $| \cdot |_{\sJ}$ of a normed ideal.

\bigskip

The dictionary entries 1. - 4. should only be viewed as a first step towards developing the noncommutative analogy as we already pointed out in \cite{27}.
One can actually go further and instead of the gradient $\nabla$ and of the rearrangement invariant norm $|| \cdot ||_{\Phi}$ consider the differential seminorm $|| \nabla(\cdot)||_{\Phi}$ to which they give rise and for which the noncommutative analogue is a noncommutative differential seminorm. Note that such a more general notion is well suited in the hybrid setting. Another remark is that the finite rank projections are just the finite projections of a type $I_{\infty}$ factor and that in the more general setting of a semifinite von Neumann algebra one should consider  finite projections of a type $II_{\infty}$ factor (\cite{27}, \cite{28}).

\bigskip

A question to which we will return later in this paper is entry 3. of the dictionary , which will require an adjustment of the analogy.

\section{The condenser quasicentral modulus}
\label{sec 13}

Applying the dictionary of the noncommutative analogy to the nonlinear condenser capacity associated with a rearrangement invariant norm produces the condenser quasicentral modulus. For the sake of simplicity our discussion will be limited to the type $I_{\infty}$ context of $\sB (\sH)$.

\bigskip

Let $\tau = (T_j)_{1 \leq j \leq n} $ be an $n$-tuple of bounded operators on $\sH$, $(\sJ, | \cdot |_{\sJ}) $ a normed ideal and $P, Q$ orthogonal finite rank projections $PQ = 0$. An application of the dictionary to the capacity quantities with respect to a rearrangement invariant norm yields the noncommutative condenser quasicentral modulus quantities:

\[
{k_{\sJ} (\tau; P, Q) = \inf \{ |[a, \tau]|_{\sJ} \vert A \in \sR^{+}_1 , AP = P, AQ = 0 \}} .
\]
\[
{k_{\sJ} (\tau ; P) = k_{\sJ} (\tau ; P , 0) } 
\]
\[
{k_{\sJ} (\tau) = \sup \{k_{\sJ} (\tau ; P) \vert P \in \sR^{+}_{1} , P^{2} = P \} }. 
\]

\bigskip

Note that the number $k_{\sJ} (\tau)$ defined above is easily seen to be equal to the quasicentral modulus discussed previously. This also justifies the use of the same notation for both.

\bigskip

Let us add to this application of the dictionary some remarks about entry 3. of the dictionary. It is a weak point of the analogy: the $\ast$-algebras $\sR$ and  $\sC^{\infty}_{0} (\Omega)$ are rather different objects. There is however a simple remedy, instead of the nonlinear potential theory on the open set $\Omega \subset \mathbb {R}^{n}$, one should consider discrete potential theory on a graph, like for instance on the Cayley graph of a group with a finite generator. Potential theory on graphs (see \cite{1} , \cite{17}) is a well-developed subject and there is also nonlinear potential theory in this setting (the appendix in \cite{17} provides a guide to some of the work in this area). With this replacement of nonlinear potential theory on $\mathbb {R}^{n} $ with the corresponding theory on suitable graphs, the noncommutative analogy becomes perfect. This will be illustrated in the next section with the case of the Cayley graphs of discrete groups with a finite generator.

\section{Discrete groups}
\label{sec 14}

Let $G$ be a discrete group with a finite generator $\gamma = (g_1, ... , g_n)$ and let $\lambda$ be the left regular representation of $G$ on $\ell^{2} (G)$. This will give rise to "classical" condensers on the Cayley graph of $G$ and to noncommutative condensers in $\sB (\ell^{2} (G))$ with respect to $\lambda (\gamma)$. Let further $(\sJ , | \cdot |_{\sJ}) $ be a normed ideal and let $\ell_{\sJ} (G) $ be the rearrangement invariant space of functions $f : G \rightarrow \mathbb{C} $ with the $\sJ$-norm i.e. the norm of the corresponding diagonal operators on $\sJ (\ell^{2} (G)) $.

\bigskip

Let $X, X_1, X_2  \subset G $ be finite subsets so that $X_1 \cap X_2 = \emptyset $. The classical nonlinear capacity quantities are :

\[
{cap_{\sJ} (X_1 , X_2) = } 
\]
\[
{ = \inf \{ \max_{1 \leq j \leq n} | u(g_j \cdot ) - u(\cdot)|_{\sJ} \vert 0 \leq u \leq 1, \text{supp u finite} , u \vert_{X_1} \equiv 1, u\vert_{X_2} \equiv 0 \} }
\]
\[
{cap_{\sJ} (X) = cap_{\sJ} (X , \emptyset)} 
\]
\[
{cap_{\sJ} (G) = \sup \{ cap_{\sJ} (X) \vert X \subset G , X  \text{ finite} \}}.
\]

\bigskip

On the other hand we have quasicentral condenser capacity quantities arising from $\lambda (\gamma)$ and the noncommutative condensers defined by projections $P_{\ell^{2} (X)}$, 
$P_{\ell^{2} (X_1) } $, $P_{\ell^{2} (X_2) }$ . We have the following result:

\[
{cap_{\sJ} (X_1, X_2) = k_{\sJ} (\lambda (\gamma) ; P_{\ell^{2}(X_1)} , P_{\ell^{2}(X_2)} ) } 
\]
\[
{cap_{\sJ} (G) = k_{\sJ} (\lambda (\gamma))}.
\]

\bigskip

Note that $cap_p (G) > 0$ and $cap_p (G) = 0$ correspond to the $p$-hyperbolicity and respectively $p$-parabolicity in the sense of Yamasaki (\cite {30}).

\section{A quasicentral condenser variational problem (\cite {28} ) }
\label{sec 15}

With the quasicentral modulus extended to condensers, next in the exploration of the analogy with nonlinear potential theory is to consider variational problems. We take a look at one of the simplest such problems associated with a condenser $(P, Q)$ an ideal $\sC_{p} , 2 \leq p < \infty$ and an $n$-tuple $\tau = (T_j)_{1 \leq j \leq n} $ of bounded selfadjoint operators $T_j = T^{\ast}_j$ . It will be convenient to consider a smoother version of the modulus by replacing the $\max_{1 \leq j \leq n} |[T_j , A] |_p $ with 
\[
\vert
\begin{pmatrix}
[T_1 , A] \\
\vdots \\
[T_n , A] 
\end{pmatrix}
\vert_{p} .
\]

Associated with the condenser are two sets 

\[
\sC^{0}_{PQ} = \{ A \in \sR^{+}_{1} \vert AP=P , AQ=0 \}  ,\\
\sC_{PQ} = \{0 \leq B \leq I \vert BP=P, BQ=0 \}  .
\]

on which $I(\cdot)$ is defined. It turns out that

\[
k_{p} (\tau) = 0 \\
\implies \inf \{ I(A) \vert A \in \sC^{0}_{PQ} \} = \inf \{ I(B) \vert B \in \sC_{PQ} \}.
\]

\bigskip

About the minimization of $I(\cdot)$  we obtained the following existence and uniqueness result .

\[
{\exists X \in \sC_{PQ}}  
\]
\[
{I(X) = \inf \{I(A) \vert A \in \sC^{0}_{PQ} \}}. 
\]

$X$ is unique mod the commutant $N$ of $\tau$.

\bigskip

If $X$ is a minimizer of $I(\cdot)$ in $\sC_{PQ}$ we found that

\[
\exists P_1 , Q_1 \text{ projections} , P \leq P_1, Q \leq Q_1 , P_{1} Q_{1} = 0, XP_1 = P_1 , XQ_1 = 0 
\]
so that 
\[
{(I - P - Q_1) \Theta (I - P - Q_1) \le 0 }
\]
and
\[
{ (I - P_1 - Q) \Theta (I - P_1 - Q) \ge 0 }
\]
where
\[
{\Theta = \sum_{1 \le k \le n}[T_k , [X, T_k] (- \sum_{1 \le j \le n} [X, T_j]^2)^{\frac p2 -1} + (- \sum_{1 \le j \le n} [X , T_j]^2)^{\frac p2 - 1} [X,T_k]] .}
\]

\bigskip

The formula for $\Theta$ can be viewed up to replacing a multiplication by a Jordan algebra product as an analogue of the $p$-Laplacian
\[ 
{div (|\nabla u|^{p-2} \nabla u) }
\]

\bigskip

Note also that for the increased condenser $(P_1, Q_1)$ we have

\[
{(I - P_1 - Q_1 ) \Theta (I - P_1 - Q_1) = 0 }
\]

\section{Semifinite work}
\label{sec 16}

We have restricted our presentation to the type $I_{\infty}$ factor $\sB (\sH) $ . While the correspondent of normed ideals, also known as symmetric spaces of operators (\cite {15}) have been thoroughly studied in the general, semifinite setting, the perturbation theory is at an earlier stage when compared to the type $I$ work and there are also some results about the quasicentral modulus. Here are a few references \cite{2} , \cite {14} , \cite {27} , \cite {28} . It  is possible that more advanced von Neumann algebra theory will need to be used in further developments.


\end{document}